\documentstyle[doublespace,leqno,11pt]{article}

    \csname @addtoreset\endcsname{figure}{section}
    \csname @addtoreset\endcsname{table}{section}

\newcounter{thanksnum}
\def\thanksnumber#1
{\setcounter{thanksnum}{\value{footnote}}\setcounter{footnote}{#1}%
                     \addtocounter{footnote}{-1}\footnotemark
                     \setcounter{footnote}{\value{thanksnum}}}
\def\newtheoremz#1{\@ifnextchar[{\@othmz{#1}}{\@nthmz{#1}}}

\def\@nthmz#1#2{%
\@ifnextchar[{\@xnthmz{#1}{#2}}{\@ynthmz{#1}{#2}}}

\def\@xnthmz#1#2[#3]{\expandafter\@ifdefinable\csname #1\endcsname
{\@definecounter{#1}\@addtoreset{#1}{#3}%
\expandafter\xdef\csname the#1\endcsname{\expandafter\noexpand
  \csname the#3\endcsname \@thmcountersepz \@thmcounterz{#1}}%
\global\@namedef{#1}{\@thmz{#1}{#2}}\global\@namedef{end#1}{\@endtheoremz}}}

\def\@ynthmz#1#2{\expandafter\@ifdefinable\csname #1\endcsname
{\@definecounter{#1}%
\expandafter\xdef\csname the#1\endcsname{\@thmcounterz{#1}}%
\global\@namedef{#1}{\@thm{#1}{#2}}\global\@namedef{end#1}{\@endtheoremz}}}

\def\@othmz#1[#2]#3{\expandafter\@ifdefinable\csname #1\endcsname
  {\global\@namedef{the#1}{\@nameuse{the#2}}%
\global\@namedef{#1}{\@thmz{#2}{#3}}%
\global\@namedef{end#1}{\@endtheoremz}}}

\def\@thmz#1#2{\refstepcounter
    {#1}\@ifnextchar[{\@ythmz{#1}{#2}}{\@xthmz{#1}{#2}}}

\def\@xthmz#1#2{\@begintheoremz{#2}{\csname the#1\endcsname}\ignorespaces}
\def\@ythmz#1#2[#3]{\@opargbegintheoremz{#2}{\csname
       the#1\endcsname}{#3}\ignorespaces}

\def\@thmcounterz#1{\noexpand\arabic{#1}}
\def\@thmcountersepz{.}
\def\@begintheoremz#1#2{ \trivlist \item[\hskip \labelsep{\bf #1\ #2}]}
\def\@opargbegintheoremz#1#2#3{ \trivlist
      \item[\hskip \labelsep{\bf #1\ #2\ (#3)}]}
\def\@endtheoremz{\endtrivlist}

\def\F{{\cal F}}

\def\R{{\bf R}}
\def\E{{\bf E}}
\def\P{{\bf P}}

\def\b{\beta}
\def\s{\delta}

\def\ww{\widetilde}

\def\s{\sigma}

\def\p{\partial}

\newcommand{\be}{\begin{equation}}
\newcommand{\ee}{\end{equation}}
\newcommand{\bd}{\begin{displaymath}}
\newcommand{\ed}{\end{displaymath}}
\newcommand{\ba}{\bd\begin{array}{rl}}
\newcommand{\ea}{\end{array}\ed}

\date{\  }
\def\s{\sigma}

\def\R{{\bf R}}

\title{
On first exit times for homogeneous diffusion processes\thanks{
Theory of Probability and Its Application, 1987, v.31, No 3, pp.
497-498; English translation of the September 1986 issue of the
Soviet journal "Teoria Veroyatnostei i ee primneniya"; translated
from Russian by Amal Ellis}}
\author{Nikolai Dokuchaev\thanks{Correspondence address:
The Institute of Mathematics and Mechanics, St.Petersburg State
University, Bibliotechnaya pl.2, Petrodvoretz, St.Petersburg,
198904,  Russia. E-mail: dokuchaev@pobox.spbu.ru; fax
7(812)4286998} }
\begin{document}
\maketitle In this paper we derive an upper bound for the
expression $\E|T_1-T_2|$, first exit times from a region for two
homogeneous diffusion processes. \par Let $(\Omega,{\cal F},\P)$
be a standard probability space and $Q$ be a bounded simply
connected region in $\R^n$. Consider the $n$- dimensional
diffusion processes $y_i(t)$, i=1,2 (see [1, v.3,pp.257-258]) of
the following It\^o equations:
\begin{eqnarray}
\label{(1)}
dy_i(t)=f_i(y(t))dt+\b_i(y(t))dw_t,\\ \label{(2)}
y_i(0)=a_i.
\end{eqnarray}
Here $t\ge 0$, $w_t$ is a standard $d$-dimensional Wiener process,
coordinated as usual with some right-continuous non-decreasing
flow of $\s$-algebras $\F_t\subset \F$; $f_i$ and $\b_i$ are
nonrandom functions with respective values in $\R^n$ and
$\R^{n\times d}$, (here and elsewhere $i=1,2$). The random vectors
$a_i$ are measurable with respect to the $\s$-algebra $\F_0$ and
$a_i\in\bar Q$ with probability 1($\bar Q$ denotes the closure of
the region $Q$ ). All vectors and matrices are real with Euclidean
norm $|\cdot|$.
\par
Consider the r.v.'s $T_i=\inf\{t:\ y_i(t)\notin Q\}$. These are
the first exit times of the processes $y_i(t)$ from the region
$Q$.
\par
Consider in $Q$ the Dirichlet problem
\be
\label{(3)} L_iv_i=-1,\quad v_i|_{\p Q}=0. \ee Here $\p Q$ is the
boundary of $Q$ and the differential operators \be\label{(4)}
L_i=\sum_{j=1}^nf_i^{(j)}\frac{\p}{\p
y^{(j)}}+\frac{1}{2}\sum_{j,k=1}^nb_i^{(j,k)}\frac{\p^2}{\p
y^{(j)}\p y^{(k)}}, \ee where $f_i^{(j)}$, $y^{(j)}$,
$b_i^{(j,k)}$ are the components of the vectors $f_i$, $y$ and the
matrices $b_i(\cdot)=\b_i(\cdot)\b_i(\cdot)^\top$.
\par
We assume that all the components of the function $f$, $\b_i$ are
continuously differentiable, the eigenvalues of the matrices $b_i$
are isolated from zero uniformly in $Q$, and the boundary $\p Q$
is smooth. Problem (3)-(4) then has  a (unique) solution that is
twice continuously differentiable up to the boundary. \par Denote
by $dv_i/dy$ the vector with components $\p v_i/\p y^{(j)}$. \par
{\bf Theorem}.{\it \be \label{(5)} \E |T_1-T_2|\le \max_{i=1,2} \,
\sup_{x\in Q}\left|\frac{d v_i}{d y}(y)\right| \E|y_1(\ww
T)-y_2(\ww T)|, \ee
 where $\ww T=T_1\land T_2=\min(T_1,T_2)$. }
 \par
 {\it Proof}.
 Let $e_1$ and $e_2$ be the indicator functions of the events
 $\{T_1>T_2\}$ and $\{T_2>T_1\}$ respectively.
 Obviously,
 \be
 \label{(6)}
\E|T_1-T_2|=\E\{e_1(T_1-T_2)\}+\E\{e_2(T_2-T_1)\}.
 \ee
 The r.v.'s $e_i$ are measurable with respect to the $\s$-algebras
 $\F_{\ww T}$, $\F_{T_j}$, associated with the Markov times (with
 respect to the flow $\F_t$) $\ww T$, $T_j$, $j=1,2$ (see [1, v.2,
 Chap.4, \S 2]). Using It\^ o's formula, we obtain the equality
 \be
 \label{(7)}
 \begin{array}{ll}
\E\{e_1\{v_1[y_1(T_2)]-v_1[y_1(T_2)]\}\}&=-
\E\{e_1\{v_1[y_1(T_1)]-v_1[y_1(T_2)]\}\}\\
&=-\E\left\{e_1\int_{\ww T}^{T_1}L_1v_1[y_1(t)]dt\right\}
=\E\{e_1(T_1-T_2)\}.
 \end{array}
 \ee
 If we replaced the indices $1,2$ in (7) by 2,1, we get an
 analogous expression for $\E\{e_1(T_1-T_2)\}$. Now the assertion of
 the theorem follows at once from (6).
 \par
 {\bf Example.} In equations (1)-(2) let $n=d=1$,
 $y_i(t)=a_i+w_t$,
 $Q=(0,1)$. Then, obviously, $v_i(y)=y(1-y)$, $dv_i(y)/dy=1-2y$ and
 $\E|T_1-T_2|\le |a_1-a_2|$.
 \par
 The behavior of the first exit times of the process when the
 drift and diffusion coefficients vary has been studied in details
 in the book [2]. However, it solves problems of different sort
 than the one considered above.
 \par
 For non-homogeneous processes and $T_1$, $T_2$ bounded times an
 estimate similar to (5)
 was derived in [3] and used to prove the "stochastic maximum
 principle" established in [4]
 for a controlled diffusion process in a bounded region. For
 strong Markov processes of arbitrary form, an estimate of type
 (5) was derived in [5] (however, the assumptions made there are
 restrictive and not easy to check for nondiffusional processes).

\end{document}